\documentclass[letterpaper, 10 pt, conference]{ieeeconf}
\IEEEoverridecommandlockouts
\usepackage{amsmath,amssymb,amsfonts}
\usepackage{graphicx}
\usepackage{booktabs}
\usepackage{xcolor}
\usepackage{cite}
\usepackage{url}

\newcommand{\R}{\mathbb{R}}

\newcommand{\Q}{\mathbb{Q}}
\newcommand{\norm}[1]{\left\lVert #1 \right\rVert}
\newcommand{\Linfnorm}[1]{\norm{#1}_{L_\infty}}
\newcommand{\sigmax}{\overline{\sigma}}
\newcommand{\conj}[1]{\widetilde{#1}}

\title{\LARGE \bf Exact Decisions and Certified Values in $L_\infty$-Norm Computation \\ from State-Space Representations}

\author{
Alban Quadrat$^{1}$, Fabrice Rouillier$^{1}$ and Grace Younes$^{2}$
\thanks{$^{1}$A. Quadrat and F. Rouillier are with Sorbonne Universit\'e and Universit\'e Paris Cit\'e, CNRS, IMJ-PRG, Inria Paris, PACE, France.
{\tt \{alban.quadrat, fabrice.rouillier\}@inria.fr}}
\thanks{$^{2}$G. Younes is with Sorbonne University Abu Dhabi, SAFIR, Abu Dhabi, UAE.
{\tt grace.younes@sorbonne.ae}}
}

\begin{document}

\maketitle
\thispagestyle{empty}
\pagestyle{empty}

\begin{abstract}
Computing the $L_\infty$-norm of a time-invariant linear ordinary differential control system is a routine task in robust control, usually delegated to floating-point routines that perform well on most practical systems. Behind the numerical value, however, lies a discrete decision: whether a pole lies on the imaginary axis, and whether a mode of a realization is actually a pole of the transfer matrix. Floating-point arithmetic can only replace this decision with a tolerance. This paper explains, for a control audience, what reliable symbolic-numeric computation offers for this problem: exact arithmetic where a decision has to be made, and certified enclosures where a value has to be computed. We show how the problem is turned into a small polynomial system directly from state-space data, how this system is solved in practice, and what it costs. A benchmark on the {\tt COMPleib} library and on families of systems with known norms shows that two standard numerical routines err in opposite directions on the imaginary-axis decision, and never on the same system, while the certified method decides it correctly on every instance tested and computes the norm at a cost comparable to that of the numerical routines on most systems. We conclude with practical guidelines for users of numerical tools.
\end{abstract}


\section{Introduction}

In robust control, the $L_\infty$-norm of a time-invariant linear system quantifies how much the system can amplify, in the worst case, the energy of a signal; for a stable system it coincides with the 
$H_\infty$-norm on which robust stability analysis, robustness margins and $H_\infty$ synthesis are built \cite{DFT90,ZDG96}. For a proper rational transfer matrix $G$ without poles on the imaginary axis,
\[
  \Linfnorm{G} = \sup_{\omega \in \R} \sigmax\bigl(G(i\omega)\bigr),
\]
where $\sigmax$ denotes the largest singular value. In practice, control engineers obtain this number from routines such as \texttt{getPeakGain} and \texttt{hinfnorm} in {\sc Matlab} or the {\sc Slicot} routine \texttt{AB13DD}. These routines rely on a Hamiltonian characterisation: for a level $\gamma$ larger than $\sigmax(G(i \, \infty))$, the norm is at least $\gamma$ exactly when an associated Hamiltonian matrix has eigenvalues on the imaginary axis, which leads to a bisection on $\gamma$ \cite{BBK89} or to a quadratically convergent iteration \cite{BB90}. For alternative methods, see  \cite{Belur, benner2000, benner2012,benner2018, BB90, Bruinsma, Gahinet, Henrion, Robel, Sima2006}. On most practical systems, they are fast and accurate.

There is, however, a question hidden behind the numerical value. The $L_\infty$-norm of a proper transfer matrix is infinite exactly when $G$ has a pole on the imaginary axis. Computing the norm therefore requires deciding whether a pole lies \emph{on} the axis or merely \emph{close} to it, and, when the system is given by a non-minimal realization, whether an imaginary eigenvalue of the state matrix is actually a pole of $G$ or cancels. These are yes-or-no questions whose answer changes discontinuously with the data. A floating-point routine can only replace them by a threshold, and a threshold admits two kinds of error: an infinite norm returned as a large finite number, and a finite norm returned as infinite. Both are observed on standard benchmark systems.

Reliable symbolic-numeric methods address such questions differently. Decisions are made in exact arithmetic, where a common root or a vanishing coefficient is detected without any tolerance, and numerical values are returned as intervals that are proved to contain the exact result. As shown in \cite{KS06}, the norm can be read on a plane algebraic curve built from $G$: it is the height of its highest real critical point, where the tangent is horizontal or two branches cross, or of a horizontal asymptote, and this point can be computed with certified computer algebra tools \cite{BQRY21,YQR25}. The price is the growth of the data: exact numbers and polynomials become larger at each operation, so the art lies in building and solving a polynomial system that remains small.

This paper is intended for control engineers. It aims to explain what these methods compute, how they are used on the $L_\infty$-norm, and what they cost, rather than to give full proofs, which can be found in the companion journal paper \cite{JACA}. Section~\ref{sec:paradigms} recalls the main paradigms of reliable computation and how symbolic and numerical methods complement each other. Section~\ref{sec:construction} shows how the polynomial system is built directly from state-space data, and Section~\ref{sec:solving} shows how it is solved in practice with the library {\tt PACE}. Section~\ref{sec:benchmark} reports a benchmark against standard numerical routines and discusses the limitations of the method, and Section~\ref{sec:conclusion} draws practical conclusions.
\section{Reliable symbolic-numeric computation}\label{sec:paradigms}

This section recalls, without proofs, the ideas on which the rest of the paper relies. More details can be found in \cite{BQRY21,YQR25,JACA}.

\subsection{Three kinds of arithmetic}

\emph{Floating-point arithmetic} approximates every real number by a number with a fixed number of binary digits ($53$ in double precision, i.e.\ about $16$ significant decimal digits). Each operation is fast and takes the same memory, but its result is rounded. The results are excellent approximations in most cases, but a floating-point routine cannot tell whether a quantity is exactly zero: a computed eigenvalue with real part $10^{-17}$ may come from a pole on the imaginary axis perturbed by rounding, or from a genuine pole close to it. Any decision of this kind must be taken with respect to a tolerance.

\emph{Exact arithmetic} works with integers, rational numbers and polynomials with rational coefficients, without any rounding. A coefficient that vanishes is exactly zero, and two polynomials that share a root have a nontrivial greatest common divisor (gcd), which is computed exactly. The price is that the size of the data grows: adding two fractions may add the number of digits of their denominators, the determinant of a $k\times k$ polynomial matrix can have degree up to $k$ times the degree of its entries, and eliminating a variable from a polynomial system may multiply the degrees. An exact computation is therefore only practical if the problem is formulated so that the data remain small.

\emph{Interval arithmetic} sits in between. Each number is enclosed in an interval whose endpoints are represented exactly, and operations are rounded outwards, so that the exact result is guaranteed to lie inside the computed interval. Combined with \emph{multiprecision}, where the number of digits is chosen as large as needed, the intervals can be made as narrow as required. Together with the certified isolation of real roots described below, this is how numerical values are certified.

\subsection{Certified isolation of real roots}

Let $P$ be a univariate polynomial with rational coefficients. \emph{Isolating} its real roots means computing disjoint intervals with rational endpoints, each containing exactly one real root, and together containing all of them. This can be done exactly: Sturm sequences count the real roots in an interval, and methods based on Descartes' rule of signs subdivide an interval until each piece is certified to contain zero or one root \cite{RZ04}. The isolation is therefore a proof, not an approximation. Each interval can then be refined, for instance by bisection, as far as required. For instance, the positive root of $x^2 - 2$ is isolated by $[1, 2]$ and refined to $[1.41, 1.42]$, and so on.

For a system of two polynomial equations in two unknowns, the same goal is reached in two steps. A variable is first eliminated, for instance by a resultant, which yields a univariate polynomial whose roots contain one coordinate of the solutions. The solutions are then described by a rational parametrisation, such as a \emph{rational univariate representation} (RUR) \cite{Rou99}: every solution is expressed as a rational function of a root of a univariate polynomial. Isolating the real roots of this polynomial and evaluating the rational functions on the isolating intervals gives certified boxes for all the real solutions of the system \cite{BLMPRS16}.

\subsection{Combining symbolic and numerical methods}

The two families of methods have complementary strengths. Exact symbolic computations are the right tool for \emph{decisions}: whether a polynomial has a root on the imaginary axis, whether two polynomials share a factor, how many real roots a polynomial has. Numerical computations are the right tool for \emph{values}, once the exact problem has been reduced to isolating the real roots of univariate polynomials. A reliable symbolic-numeric method uses exact arithmetic for the first and certified numerics (root isolation and interval refinement) for the second. It returns a result together with a guarantee: here, an interval that is proved to contain the $L_\infty$-norm, or a proof that the norm is infinite.

The difficulty is to keep the symbolic part small. Most of the work described in Sections~\ref{sec:construction} and~\ref{sec:solving} consists in formulating the problem so that the polynomials handed to the solver have low degree and small coefficients. The implementation used in this paper is written in {\sc Julia} and distributed within the {\tt PACE} project\footnote{\url{https://pace.gitlabpages.inria.fr/tutorials/}}.

\subsection{The $L_\infty$-norm as a problem on a plane curve}

Let $\conj{G}(s) := G^T(-s)$. The following characterisation is due to Kanno and Smith \cite{KS06}: for $\gamma > 0$, we have $\gamma > \Linfnorm{G}$ if and only if $\gamma > \sigmax\bigl(G(i\infty)\bigr)$ and
\[
  \det\bigl(\gamma^2 I - \conj{G}(i\omega)\,G(i\omega)\bigr) \neq 0 \quad \text{for all } \omega \in \R .
\]
This determinant is a rational function of $\omega$ and $\gamma$. Let $n(\omega, \gamma)$ be its numerator; it is a polynomial in $\omega^2$ and $\gamma^2$. Since $G$ has no pole on the imaginary axis, for a fixed $\omega$ the real roots $\gamma$ of $n(\omega,\gamma)$ are the $\pm\sigma_i(G(i\omega))$, so that $\Linfnorm{G}$ is the supremum of the $\gamma$-coordinates of the real points of the curve $n(\omega,\gamma) = 0$. This supremum is either attained at a real critical point of the curve, i.e.\ a real solution of
\[
  n(\omega,\gamma) = 0, \qquad \frac{\partial n}{\partial \omega}(\omega,\gamma) = 0,
\]
where $n$ is replaced by its square-free part, or approached as $\omega \to \infty$, along a horizontal asymptote given by a real root of the leading coefficient of $n$ in $\omega$ \cite{KS06,BQRY21}. The problem is thus reduced to a bivariate polynomial system, which is solved with the tools above.

\paragraph*{Example} For $G(s) = 1/(s^2 + s + 1)$, we have $|G(i \, \omega)|^2 = 1/(\omega^4 - \omega^2 + 1)$, so that
\[
  n(\omega, \gamma) = \gamma^2 \, (\omega^4 - \omega^2 + 1) - 1,
  \;
  \frac{\partial n}{\partial\omega} = 2 \, \gamma^2 \, \omega \, (2 \, \omega^2 - 1).
\]
The critical points are at $\omega = 0$, where $\gamma = 1$, and at $\omega^2 = 1/2$, where $\gamma^2 = 4/3$; the leading coefficient $\gamma^2$ gives the horizontal asymptote $\gamma = 0$. Hence $\Linfnorm{G} = 2/\sqrt{3}$, attained at $\omega = 1/\sqrt{2}$. A certified method returns this value as an exact algebraic number, or as an interval with rational endpoints of any desired width.
\section{Building the polynomial system}\label{sec:construction}

In control applications, a linear system is usually given by a state-space representation
\[
  \dot{x} = A \, x + B \, u, \qquad y = C \, x + D \, u,
\]
with $A \in \R^{n_x\times n_x}$, $B \in \R^{n_x\times m}$, $C \in \R^{p\times n_x}$, $D \in \R^{p\times m}$, and the transfer matrix $G(s) = C \, (s \, I_{n_x} - A)^{-1} \, B + D$. This section shows how the curve of Section~\ref{sec:paradigms} is obtained from these data while keeping its degrees and the size of its coefficients as small as possible. In exact arithmetic, this is where most of the effort goes: a careless construction can cost more than solving the resulting system.

\subsection{From the state-space representation, without forming the transfer matrix}

The natural route would be to invert $s \, I_{n_x} - A$ symbolically, form the $p \, m$ rational function entries of $G$, and put them over a common denominator. In exact arithmetic, this is expensive, because rational functions in $s$ grow quickly. Instead, we write $G = N/d$ with
$$\begin{aligned}
  d(s) & = \det(s \, I_{n_x} - A), \vspace{1mm} \\ N(s) & = C\,\mathrm{adj}(s \, I_{n_x} - A)\,B + d(s)\,D,
\end{aligned}$$
where $\mathrm{adj}(M)$ denotes the adjugate matrix of $M$, and compute $d$ and the polynomial matrix $N$ together by the Faddeev--Leverrier recurrence, which expresses the adjugate of $s \, I_{n_x} - A$ as $\sum_k s^k M_k$ with $M_{k-1} = A \, M_k + d_k \, I_{n_x}$. The matrices $M_k$ are never formed: the recurrence is applied to $B$ directly, and $C$ is applied on the fly, so that every intermediate object has the size of $B$ or of $C \, B$. Only polynomial arithmetic with rational coefficients is involved. If $A$ has rational entries, scaling $A$ by the least common multiple $q$ of the denominators of its entries and substituting $u = q \, s$ makes the recurrence run over the integers. On the $114$ {\tt COMPleib} models on which both routes were timed, this route is faster than forming $G$ explicitly on $113$, by a median factor of $13$ that grows with the size of the system, up to $1$\,s against $3980$\,s for a model with $120$ states \cite{JACA}.

If the user already holds $G$ as a matrix of rational functions, the pipeline starts one step later, with $d$ the least common denominator of the entries of $G$ and $N = d \, G$.

\subsection{Choosing the smaller side}

On the imaginary axis, $\conj{G}(i \, \omega) = G(i \, \omega)^\ast$, so the squared singular values of $G(i \, \omega)$ are the eigenvalues of the $m\times m$ matrix $\conj{G} \, G$, and also, apart from zeros, of the $p\times p$ matrix $G \, \conj{G}$. The curve can be built from either one; building it from the smaller one makes its degree in $\gamma^2$ equal to $\min(m,p)$ rather than $\max(m,p)$. With $e(s) = d(s) \, d(-s)$, the curve is obtained as the characteristic polynomial of the polynomial matrix $\conj{N} \, N$ (or $N \, \conj{N}$), evaluated at $\gamma^2 \, e(s)$ and read at $s = i \, \omega$: no division by a rational function is ever performed.

\subsection{Halving the degrees}

The polynomial $n(\omega,\gamma)$ only involves $\omega^2$ and $\gamma^2$. We therefore work in the variables $\Omega = \omega^2$ and $\Gamma = \gamma^2$, with $n(\omega,\gamma) = c(\Omega,\Gamma)$; in the construction above, this amounts to substituting $s^2 = -\Omega$ and $\gamma^2 = \Gamma$. This halves the degree in each variable and divides the number of solutions of the system, counted with multiplicities, by four, which reduces the cost of the solving step accordingly.

This change of variables has a subtle consequence. The frequency now ranges over $\Omega \ge 0$, and the critical points of the curve at $\omega = 0$ disappear from the polynomial system 
$$\{c = 0,\ \partial c/\partial\Omega = 0\}.$$ Yet the norm is often attained at $\omega = 0$: on the $116$ {\tt COMPleib} models of moderate size, it is attained at $\omega = 0$ on $48$, at an interior critical point on $43$, and is infinite on the remaining $25$; the boundary $\omega = 0$ is thus the most frequent case among the systems with a finite norm \cite{JACA}. The candidates for $\Linfnorm{G}^2$ are therefore three families of values of $\Gamma$:
\begin{enumerate}
  \item the \emph{interior candidates}: real solutions of the polynomial system $c = \partial c/\partial\Omega = 0$ with $\Omega \ge 0$;
  \item the \emph{boundary candidates}: real roots of the univariate polynomial $c(0,\Gamma)$, for $\omega = 0$;
  \item the \emph{asymptotic candidates}: real roots of the leading coefficient of $c$ in $\Omega$, which are the squared singular values of $D = G(i \, \infty)$, for $\omega \to \infty$.
\end{enumerate}
The norm is the largest of all these candidates \cite{JACA}. The last two families only require univariate polynomials and cost almost nothing. They are not a formality: for the closed-loop projector of the unstable plant $k/s$ stabilised by the static controller $-1/k$, whose inverse norm is the robust stability radius of the closed-loop system, the norm is $\sqrt{k^2+1}/k$ for $0<k<1$, a value attained only at infinite frequency \cite{JACA}.

\subsection{Deciding whether a pole lies on the imaginary axis}

Before building the curve, one must know whether $G$ has a pole on the imaginary axis, in which case $\Linfnorm{G} = \infty$. A realization may carry uncontrollable or unobservable modes: eigenvalues of $A$ that cancel in every entry of $N/d$ and are not poles of $G$. On {\tt COMPleib}, this happens on $37$ of the $114$ models examined; in one of them, AC14, $37$ of the $40$ modes cancel \cite{JACA}. Taking an eigenvalue of $A$ on the axis for a pole of $G$ returns an infinite norm for a system whose norm is finite.

In exact arithmetic, both questions are answered without any tolerance:
\begin{itemize}
  \item the poles of $G$ are the roots of the polynomial $d_{\min} = d / \gcd(d, N_{11}, \dots, N_{pm})$, where the $N_{ij}$'s are the entries of $N$, which removes every cancelling mode;
  \item writing $d_{\min}(i \, y) = d_{\mathrm{re}}(y) + i\,d_{\mathrm{im}}(y)$ with real polynomials $d_{\mathrm{re}}$, $d_{\mathrm{im}}$, the imaginary roots $i \, y$ of $d_{\min}$ are the common real roots of $d_{\mathrm{re}}$ and $d_{\mathrm{im}}$, hence the real roots of their gcd, which are counted exactly.
\end{itemize}
The answer is an integer, the number of poles on the axis, obtained without any threshold. Among the $116$ {\tt COMPleib} models of moderate size, this settles $25$ outright, with an infinite norm established rather than estimated. Asymptotic stability, needed to compare with routines that return $\infty$ for every system they consider unstable, is decided in the same spirit. A sign test on the coefficients of $d_{\min}$ rejects most unstable systems. The others go through a certified numerical step: the roots of $d_{\min}$ are enclosed in complex balls at increasing precision ($64$ to $512$ bits), and the verdict is read on the signs of their real parts as soon as every ball lies on one side of the axis, so that the verdict is certified (a pole on the axis has been excluded before). A Routh table in integer arithmetic settles the few remaining cases. On the benchmark, the sign test decides $38$ systems, the balls $184$ and the Routh table $5$.
\section{Solving the system with PACE}\label{sec:solving}

Once the curve is built, the bivariate system $$\{c = 0,\ \partial c/\partial\Omega = 0\}$$ is reduced by removing the common factors of $c$ and $\partial c/\partial\Omega$; this loses none of the critical points that matter \cite{JACA}. It is then solved in {\tt PACE} with exact computations followed by certified root isolation. Three choices make this step efficient.

\subsection{Which variable to eliminate}

Either $\Omega$ or $\Gamma$ can be eliminated. We eliminate $\Gamma$. The cost of the subresultant computations that perform the elimination is driven by the degree in the eliminated variable, and the degree in $\Gamma$ is only $\min(m,p)$: it equals one on $125$ of the $272$ systems of the benchmark and never exceeds $12$ on the {\tt COMPleib} models but one, whereas the degree in $\Omega$ grows with the state dimension, up to $480$. On CM4, a model with $240$ states and one input, eliminating $\Omega$ instead would multiply the cost of the elimination by $(d_\Omega/d_\Gamma)^2$, about $5.8\cdot10^4$ \cite{JACA}. Eliminating $\Gamma$ also avoids testing the roots of an eliminant in $\Gamma$ one by one for a real frequency $\Omega\ge0$ above them. The resulting univariate polynomial lives in $\Omega$, and each of its real roots gives $\Gamma$ through a rational function. As a side benefit, the frequency at which the norm is attained is returned with the same guarantee as the norm itself. The elimination itself is then nearly free: on $66$ models of the benchmark, computing the eliminant alone takes $0.06$\,s in total, against $1.01$\,s for the whole parametrisation \cite{JACA}. The cost of the solving step lies in the parametrisation, and this is why the next choice matters.

\subsection{Parametrisations without a black-box solver}

A general-purpose solver for polynomial systems, for instance one based on Gr\"obner bases, can be used as a black box, but it ignores the structure of the problem. Here the parametrisation is obtained at three levels, each one being used only when the previous one does not apply.
\begin{itemize}
  \item When the curve has degree one in $\Gamma$, i.e., when we have $c = a_1(\Omega)\,\Gamma + a_0(\Omega)$, the critical points are the roots of $a_1 \, a_0' - a_0 \, a_1'$ and $\Gamma = -a_0/a_1$ at each of them, so that no solver is needed. This is the case for $106$ of the $272$ systems of the benchmark.
  \item Otherwise, the subresultant sequence of $c$ and $\partial c/\partial\Omega$ with respect to $\Gamma$ gives a triangular decomposition of the system \cite{BLMPRS16}, and on each component $\Gamma$ is read as a rational function of $\Omega$ from the subresultant of degree one. This covers the $114$ other systems on which an interior system has to be solved.
  \item In general, a rational univariate representation (RUR) \cite{Rou99,BLMPRS16} may require a linear change of variables, i.e.,  $T = \Omega + a\,\Gamma$, which is the most expensive step of the method. It is kept as a last resort, and it was never needed: on all the $272$ systems, the coefficient $a$ is zero.
\end{itemize}
The subresultants, and the RUR when it is needed, are computed by the {\tt LACE} library modulo several primes and reconstructed over $\Q$ by the Chinese remainder theorem, which avoids the growth of the coefficients of a direct computation \cite{JACA}. Compared with our earlier pipeline \cite{YQR25}, which computed the RUR from a lexicographic Gr\"obner basis, the largest $5\times5$ instance of its random benchmark goes from $3983$\,s to $16$\,s, with the same value \cite{JACA}.

\subsection{Isolation and refinement}

Only univariate polynomials with integer coefficients are then handled numerically. All the real roots of the eliminant in $\Omega$ are isolated by a certified method, and $\Gamma$ is evaluated on each isolating interval through the parametrisation. Isolating all the real roots costs no more, in the worst case, than isolating one \cite{JACA}. The interior candidates obtained in this way are compared with the boundary and asymptotic candidates, and the largest of them is the norm. The output is an interval with rational endpoints that is proved to contain $\Linfnorm{G}^2$, together with an interval for $\Omega^\star=\omega^{\star2}$, the frequency where it is attained, which is exact when $\omega^\star = 0$.

In practice, the cost of this step is kept low by refining as little as possible, in the following way. All the real roots of each eliminant are isolated once, at a very low precision (a few bits), and $\Gamma$ and $\Omega$ are evaluated on each isolating interval by interval arithmetic through the parametrisation. Candidates with $\Omega \le 0$ are discarded by a sign test only, an interval overlapping $0$ being refined until its sign is known. The candidates are then compared on their upper bounds, and the precision is doubled only when the leading interval overlaps another one; in that case, only the leader and the intervals overlapping it are refined, the others being left as they are. Note that roots where the parametrisation reads $0/0$ are counted and that the rounds continue until none is left, so that no candidate is lost. Once the leader is separated, it is refined alone, from its isolating interval, to the requested width ($2^{-64}$ relative by default). When two candidates are equal, the boundary $\omega = 0$ is preferred, then $\omega \to \infty$, so that the frequency is returned exactly whenever it can be.

\vspace{1mm}

In {\tt PACE}, the whole pipeline is a single call of the \texttt{PaceControl} package:
{\footnotesize\begin{verbatim}
M = read_matrices("AC1_matrices.txt")
r = hinf_peak(M)     # r.lo <= gamma^2 <= r.hi
r.infinite           # a pole on the axis
r.src, r.olo, r.ohi  # where attained, omega^2
\end{verbatim}}

The computations of this paper can be reproduced with the {\tt PaceControl} tutorials available at \url{https://pace.gitlabpages.inria.fr/tutorials/}.
\section{Performance, limitations and lessons for users of numerical tools}\label{sec:benchmark}

\subsection{Setting}

The method was run on the $168$ plant models of the {\tt COMPleib} library \cite{Lei04} (release r1\_1), a standard collection of named models used in robust control, and on families of systems constructed so that their norm is known exactly, $272$ systems in total. It was compared with the {\sc Matlab R2025a} routines \texttt{getPeakGain}, called with a relative tolerance of $10^{-10}$ instead of its default $10^{-2}$, and \texttt{hinfnorm}, and with the {\tt Slicot} routine \texttt{AB13DD} called through {\tt python-control 0.10.2}. Each realization is exported once as a decimal file, read by both sides; our implementation converts each decimal number to the exact rational it denotes. All timings are wall-clock, on an Apple M4 Max laptop \cite{JACA}.

\subsection{The imaginary-axis decision}

The main finding concerns the yes-or-no question of Section~\ref{sec:construction}. The two numerical implementations fail on it in opposite directions, and on different systems (Table~\ref{tab:axis}):
\begin{itemize}
  \item \texttt{getPeakGain} returns a large finite number for $17$ systems whose norm is infinite, up to $1.45\cdot10^{17}$ on a {\tt COMPleib} model;
  \item {\tt python-control} returns $\infty$ for $13$ systems whose norm is finite. On the eight {\tt COMPleib} models among them, the realization has imaginary-axis eigenvalues that are neither controllable nor observable and cancel in $G$; {\sc Matlab} agrees with the exact value to fifteen digits on these systems.
\end{itemize}
{The remaining five failures of} {\tt python-control} occur on constructed systems with poles at distance $2^{-28}$ or less from the axis. Over the $272$ systems, the two sets of failures are disjoint. Both behaviours are the expected consequences of a threshold: one routine declares poles lying exactly on the axis to be off it, the other declares modes that cancel, or poles very close to the axis, to be poles on the axis. No single threshold avoids both errors, because the underlying question is not a matter of accuracy. The exact method decides it on every instance.

\begin{table}[t]
\centering
\caption{Errors on the imaginary-axis decision ($272$ systems)}
\label{tab:axis}
\begin{tabular}{lcc}
\toprule
 & infinite norm & finite norm \\
 & returned finite & returned $\infty$ \\
\midrule
\texttt{getPeakGain} & 17 & 0 \\
python-control (\texttt{AB13DD}) & 0 & 13 \\
certified method & 0 & 0 \\
\bottomrule
\end{tabular}
\end{table}

\subsection{Accuracy: the state basis matters, not the peak}

On the $90$ {\tt COMPleib} models with a finite norm that {\sc Matlab} also returns as finite, \texttt{getPeakGain} agrees with the certified value to better than $10^{-10}$ on $88$, with a median relative gap of $6.1\cdot10^{-16}$. The constructed families show where accuracy is lost. A sharp resonance peak, or poles very close to the axis, did not defeat the routine on a well-scaled realization. What did is the conditioning of the state basis:
\begin{itemize}
  \item the same transfer function written in an increasingly ill-conditioned basis yields errors of up to $19\%$, and standard preprocessing (\texttt{prescale}, \texttt{minreal}, \texttt{balreal}) does not repair it;
  \item a cascade of $k$ identical second-order systems, i.e., $G(s) = \bigl(c_0/(s^2 + 2 \, a \, s + c_0)\bigr)^k$, has a norm known in closed form, $(c_0/(2 \, a \, b))^k$ with $c_0 = a^2 + b^2$, yet every minimal realization carries Jordan blocks of size $k$, so no change of basis can make it well conditioned. For $k = 19$, the two numerical routines err by $2$--$3\%$, each one being the more accurate on one of two variants of the family (with a constant or a generic numerator), so that a user holding both cannot tell which to believe.
\end{itemize}
The certified method returns enclosures of relative width below $10^{-24}$ on these families \cite{JACA}.

A further remark is useful in practice. A value that is too low can be refuted by anyone, by evaluating $\sigmax(G(i \, \omega))$ at a frequency where it is larger. A value that is too high cannot be refuted; refuting it requires an upper bound valid for all frequencies, which no floating-point routine provides, and a certified enclosure does.

\subsection{Cost}

On the $116$ {\tt COMPleib} models of moderate size, all of which are settled, the median running time is $3$\,ms, $95$ are solved in less than $0.1$\,s and $107$ in less than one second. Compared with \texttt{getPeakGain}, the certified method is at parity: the median ratio of the running times is $1.05$. The running time is governed by the product $\min(m,p)\,n_x$, which is known before any computation. The largest system solved has $480$ states ($954$\,s). On the CM family, where only $n_x$ varies, from $20$ to $480$, the running time grows like a power of $n_x$ whose local exponent rises from $2.8$ to $4.65$, against the $n_x^4$ predicted by the bounds on the sizes of the integers; at $n_x = 480$, the curve carries integers of $42\,755$ decimal digits \cite{JACA}. Building the curve accounts for $57.5\%$ of the total time, the RUR and the root isolation for $38.6\%$, and the stability decisions for $3.9\%$ \cite{JACA}.

\subsection{Limitations}

The certified method has its own limits. Its cost grows with the product $\min(m,p)\,n_x$, and on large systems, building the curve becomes the dominant step: exact polynomials and coefficients grow with the state dimension, even with the choices of Section~\ref{sec:construction}. The method works with exact rational data: a model given with floating-point coefficients is converted exactly, so that the result is certified for the model as given, not for the physical system it approximates. When the data are given in decimal form, the common denominator is a power of ten, and the integers can reach thousands of digits on systems with a few hundred states, which slows down the exact computations without affecting the result. This construction is also the one step to which a multimodular scheme has not yet been applied. Beyond a certain size, the curve itself becomes too large to be written down: on {\tt COMPleib}, four models with up to $960$ states are not settled, and $41$ further square models, mostly discretised two-dimensional heat equations with $\min(m,p)=n_x$, are out of reach \cite{JACA}. Finally, the method computes the norm for numerical values of the coefficients; when the system depends on parameters, describing the norm as a function of them requires other tools \cite{QRY24}, which are the subject of ongoing work.

\subsection{Guidelines}

For users of numerical tools, the benchmark suggests the following.
\begin{enumerate}
  \item On a well-conditioned realization of a system whose poles are clearly away from the imaginary axis, standard routines are accurate and fast.
  \item A very large finite value may be an infinite norm in disguise. When the model contains integrators, undamped modes, or structural symmetries, the presence of poles on the axis should be decided exactly.
  \item An infinite value returned for a non-minimal realization may be wrong: imaginary-axis eigenvalues of $A$ may cancel in $G$. Removing them with \texttt{minreal} again relies on a tolerance.
  \item Rescaling does not recover accuracy lost to an ill-conditioned state basis.
  \item When two routines disagree, the lower value can be checked by evaluating $\sigmax(G(i \, \omega))$; the higher one can only be confirmed by a certified upper bound.
\end{enumerate}

\section{Conclusion}\label{sec:conclusion}

Computing the $L_\infty$-norm involves a decision, whether a pole lies on the imaginary axis and whether a mode is a pole at all, that floating-point arithmetic can only settle with a tolerance. Reliable symbolic-numeric methods settle it exactly with a few gcds and real-root counts, and return the norm as a certified enclosure. The art lies in formulating the problem so that the exact computations stay small: building the curve from the state-space data, using the smaller side, halving the degrees while keeping the boundary candidates, and eliminating the right variable. With these choices, the certified computation is at parity with standard numerical routines at the median of the benchmark. Numerical methods remain far more scalable, but when the answer depends on a yes-or-no decision rather than on a magnitude, a tolerance cannot replace an exact answer.

The approach extends to descriptor systems, where the transfer matrix may be improper, and to discrete-time systems, through a M\"obius transform that maps the unit circle onto the imaginary axis \cite{JACA}. Physical models also depend on parameters, and describing the norm as a function of them \cite{QRY24} is the subject of ongoing work.

\end{document}